\documentclass[10pt]{amsart}
\usepackage{amssymb,amsfonts,amsmath}
\usepackage[all]{xy}  %for xypic 

\title{Zero cycles on generic hypersurfaces of large degree}

\author{Najmuddin Fakhruddin}

\address{School of Mathematics\\Tata Institute of Fundamental research\\Homi Bhabha Road\\
Mumbai 400005\\India}

\email{naf@math.tifr.res.in}

\newcommand{\nc}{\newcommand}

\nc{\ti}{\times}

\nc{\noi}{\noindent}
\nc{\Q}{{\bf{\Bbb{Q}}}}
\nc{\R}{\Bbb{R}}
\nc{\C}{\Bbb{C}}
\nc{\bP}{\Bbb{P}}
\nc{\Z}{\Bbb{Z}}
\nc{\V}{\Bbb{V}}
\nc{\W}{\Bbb{W}}
\nc{\Ha}{\Bbb{H}}
\nc{\F}{\Bbb{F}}
\nc{\G}{\Bbb{G}}
\nc{\gog}{\frak{g}}
\nc{\gon}{\frak{n}}
\nc{\gom}{\frak{m}}
\nc{\gop}{\frak{p}}
\nc{\got}{\frak{t}}

\nc{\cp}{\rm C^{\prime}}
\nc{\et}{\tilde E}
\nc{\tp}{\tilde p}

\nc{\la}{\lambda}
\nc{\de}{\delta}
\nc{\De}{\Delta}
\nc{\Om}{\Omega}
\nc{\Ga}{\Gamma}
\nc{\ga}{\gamma}
\nc{\si}{\sigma}
\nc{\na}{\nabla}
\nc{\tth}{\Theta}

\nc{\rar}{\rightarrow}
\nc{\Rar}{\Rightarrow}
\nc{\lrar}{\longrightarrow}

\nc{\cHc}{\mathcal{H}}
\nc{\K}{\mathcal{K}}
\nc{\lc}{\mathcal{L}}
\nc{\Cc}{\mathcal{C}}
\nc{\Od}{\mathcal{O}_D}
\nc{\Vc}{\mathcal{V}}
\nc{\Ac}{\mathcal{A}}
\nc{\I}{\mathcal{I}}
\nc{\X}{\mathcal{X}}
\nc{\ds}{\bigoplus}
\nc{\te}{\otimes}
\nc{\inters}{\bigcap}
\nc{\nab}{\bigtriangledown}

\nc{\om}{\Omega}
\nc{\omdone}{\omd^1}
\nc{\omdtwo}{\omd^2}

\nc{\chdx}{\rm CH^d (X)}
\nc{\chdxq}{\rm CH^d (X)_{\Q}}
\nc{\chdox}{\rm CH^d_0 (X)}
\nc{\chdsx}{\rm CH^d_s (X)}
\nc{\adx}{\rm A^d (X)}
\nc{\abdx}{\rm Ab^d (X)}
\nc{\htdq}{\rm H^{2d}_{et} (X, \Q_l)}
\nc{\htq}{\rm H^2_{et} (X, \Q_l)}
\nc{\bx}{\rm B^1 (X)}
\nc{\bdx}{\rm B^d (X)}
\nc{\bdox}{\rm B^{d+1} (X)}

\newtheorem{prop}{Proposition}
\newtheorem{lem}{Lemma}
\newtheorem{thm}{Theorem}

\newtheorem{cor}{Corollary}

\theoremstyle{remark}

\begin{document}

\input xypic

\begin{abstract}
We show that given a smooth projective variety $X$ over $\mathbf{C}$
with $\dim(X) \geq 3$,
 an ample line bundle $\mathcal{O}(1)$ on $X$
and an integer $n > 1$, any $n$ distinct points on a generic hypersurface of degree $d$
are linearly independent in $CH_0(X)$ if $d >> 0$. 
This generalizes a result
of C.~Voisin.
\end{abstract}
\maketitle

Let $X$ be a smooth projective algebraic variety over $\C$ of dimension
$r+1 \geq 3$. Let ${\mathcal{L}}$ be an
ample line bundle on $X$. Let $S^d = H^0 (X,{\mathcal{L}}^d)$ and for a point $x \in X$,
let $S^d_x = H^0(X,{\mathcal{L}}^d \te {\mathcal{I}}_x)$, where ${\mathcal{I}}_x$ is the ideal sheaf of $x$.
 Let $R = \oplus_{d \geq 0} S^d $.

\begin{lem}
There exists an integer $m > 0$ such that, for all $d>>0$, the following holds:

\noi (i) The natural map $S^d \te S^m \rar S^{d+m}$ is surjective.

\noi (ii) The natural map $S^d \te S^m_x \rar S^{d+m}_x$ is surjective.
\end{lem}
\begin{proof}
Let the ring $R$ be generated in degrees $\leq t$. Let $m_1 = t!$ . It is easy to 
see that $a\cdot m_1$ satisfies condition (i) of the lemma for $d> a\cdot t!\cdot t(t+1)/2$, for any positive integer $a$.
Let $d'$ be a positive integer so that  ${\mathcal{L}}^{d'}$ is generated
by its global sections.
There exists $m$, a positive multiple of $m_1$, such that the maps
$S^{d'} \te S^{m}_x \rar S^{d' + m}_x$ are surjective for all $x \in X$. 
For $d>> 0$, both the maps $S^{d-d'} \te S^{d'} \rar S^d$ and
$S^{d-d'} \te S^{d' + m}_x \rar S^{d + m}_x$ are  surjective since ${\mathcal{L}}$ is ample, hence the map
$S^{d-d'} \te S^{d'} \te S^{m}_x \rar  S^{d + m}_x$ is also surjective.
The following commutative diagram then shows that the map $S^d \te S^{m}_x \rar S^{d + m}_x$ is surjective.
$$
\diagram
S^{d-d'} \te S^{d'} \te S^{m}_x  \dto \rto  & S^{d-d'} \te S^{d' + m}_x \dto \\
S^d \te S^{m}_x \rto  & S^{d + m}_x \\
\enddiagram
$$
\end{proof}

Let ${\mathcal{X}} \subset X \times S^d$ be the universal hyperplane section. For $s \in
S^d$ we denote the fibre $p_2^*(s)$ by $X_s$, which we shall assume to be
smooth. For a vector bundle ${\mathcal{V}}$ on
$X$, by ${\mathcal{V}} (b)$ we shall mean ${\mathcal{V}} \te {\mathcal{L}}^b$.

\begin{prop}
For $d >> 0$, the bundle $T{\mathcal{X}}(m)_{|X_s}$ is generated by its global sections.
\end{prop}
\begin{proof}
The proof, given the previous lemma, is  identical to Proposition 1.1 of \cite{voisin} and is hence omitted.
\end{proof}

\begin{cor}
There exists a linear function $d(n)$ of n, such that for all $d \geq d(n)$,
the vector bundle $\om^{dim S_d}_{{\mathcal{X}}}|_{X_s}$ separates any $n$ distinct 
points of $X_s$ i.e. the global sections of the bundle surject onto the 
global sections of the bundle restricted to any subscheme consisting of $n$
distinct reduced points.
\end{cor}
\begin{proof}
 $\om^{dim S_d}_{{\mathcal{X}}}|_{X_s} \cong {\om^r_{{\mathcal{X}}}|_{X_s}}^* \te {{\mathcal{K}}_{{\mathcal{X}}}}|_{X_s} \cong \wedge^r T{\mathcal{X}} |_{X_s} \te {\mathcal{K}}_{X_s} \cong  \wedge^r T{\mathcal{X}}(r\cdot m)|_{X_s} \te {\mathcal{K}}_{X_s}(-r \cdot m) \cong  \wedge^r T{\mathcal{X}}(r\cdot m) |_{X_s} \te {\mathcal{K}}_X(d -r \cdot m)|_{X_s}$. By the proposition $ \wedge^r T{\mathcal{X}}(r\cdot m) |_{X_s}$ is generated
by global sections if $d>>0$. Since ${\mathcal{L}}$ is ample, there exists a linear function
$d(n)$ such that $ {\mathcal{K}}_X(d -r \cdot m) |_{X_s}$ separates $n$ points if $d \geq d(n)$. It follows that the tensor product also separates $n$ distinct points.
\end{proof}

\begin{thm}
Let $X$ be a smooth projective variety of dimension $r+1 \geq 3$ and let ${\mathcal{L}}$ be an ample line bundle on $X$. Then there exists a linear function $d(n)$
such that for all $d \geq d(n)$, any $n$ distinct points of a generic 
hypersurface $X_s$, $s \in S^d$, are linearly independent in $CH^r (X_s)$.
\end{thm}
\begin{proof}
Suppose not. Then there exists an etale map $S \rar S^d$ and $n$ distinct
sections $\si_1,...,\si_n$ of ${\mathcal{X}}_S$ such that the classes of these sections in
$CH^r ({\mathcal{X}}_S)$ are linearly dependent. We may assume that $S$ is affine and that all the fibres are smooth. Consider the classes of these cycles,
$[\si_i]$, in the
Hodge cohomology group $H^r ({\mathcal{X}}_S, \om^r_{X_S})$. By the Grothendieck-Serre duality, it is easy to see that as an element of $Hom(H^0 ({\mathcal{X}}_S,\om^{dim S}_{{\mathcal{X}}_S}),H^0 (S, \om^{dim S}_S))$, $[\si_i]$ is nothing but the restriction map $\si_i^*$ on
differential forms. By the previous corollary we see that all the  $\si_i^*$
are linearly independent, which is a contradiction.
\end{proof}

\begin{cor}
Let $X$, ${\mathcal{L}}$ be as above. There exists a linear function $d'(n)$ such that
for all $d \geq d'(n)$, the generic hypersurface $X_s$, $ s \in S^d$, does not
contain any (possibly singular) $n$-gonal curves.
\end{cor}
\begin{proof}
Follows easily from the theorem by considering two distinct elements in the
linear system corresponding to a degree $n$ map from the normalisation of
the curve to ${\bP}^1$.
\end{proof}

\emph{Acknowledgements.}
This note was written while the author was partially supported by an N.S.F.~grant
at the Institute for Advanced Study, Princeton, during the year 1995-1996.

\end{document}